\documentclass[11pt]{article}
\usepackage[cm]{fullpage}

\usepackage{color}

\usepackage{amsmath,amscd, float,rotating}
\usepackage{pb-diagram}

\usepackage{latexsym}
\usepackage{amsfonts}
\usepackage{amsmath}
\usepackage{amssymb}

\numberwithin{equation}{section}
\newcommand{\qed}{\hfill \ensuremath{\Box}}

\def\XXint#1#2#3{{\setbox0=\hbox{$#1{#2#3}{\int}$}
\vcenter{\hbox{$#2#3$}}\kern-.5\wd0}}

\newcommand{\dbar}{\overline{\partial}}

\newcommand{\ddt}[1]{\frac{\partial #1}{\partial t}}

\newcommand{\ov}[1]{\overline{#1}}

\newcommand{\ddbar}{\frac{\sqrt{-1}}{2\pi} \partial\dbar}

\begin{document}
\newcounter{remark}
\newcounter{theor}
\setcounter{remark}{0} \setcounter{theor}{1}
\newtheorem{claim}{Claim}
\newtheorem{theorem}{Theorem}[section]
\newtheorem{proposition}{Proposition}[section]
\newtheorem{lemma}{Lemma}[section]
\newtheorem{definition}{Definition}[section]
\newtheorem{conjecture}{Conjecture}[section]
\newtheorem{corollary}{Corollary}[section]
\newenvironment{proof}[1][Proof]{\begin{trivlist}
\item[\hskip \labelsep {\bfseries #1}]}{\end{trivlist}}
\newenvironment{remark}[1][Remark]{\addtocounter{remark}{1} \begin{trivlist}
\item[\hskip \labelsep {\bfseries #1
\thesection.\theremark}]}{\end{trivlist}}
\newenvironment{example}[1][Example]{\addtocounter{remark}{1} \begin{trivlist}
\item[\hskip \labelsep {\bfseries #1
\thesection.\theremark}]}{\end{trivlist}}
~

\centerline{\bf \LARGE Metric flips with Calabi ansatz \footnote{The first named author is supported in part by National Science Foundation grant DMS-0847524 and a Sloan Foundation Fellowship.}}

\bigskip

\begin{center}{\large Jian Song$^*$ ~  and  ~  Yuan Yuan$^\dagger$ }
\end{center}

\bigskip

\noindent
{\bf Abstract} We study the limiting behavior of the K\"ahler-Ricci flow on $\mathbb{P}(\mathcal{O}_{\mathbb{P}^n} \oplus \mathcal{O}_{\mathbb{P}^n}(-1)^{\oplus (m+1)})$ for $m, n\geq 1$, assuming the initial metric satisfies the Calabi symmetry. We show that the flow either shrinks to a point, collapses to $\mathbb{P}^n$ or contracts a subvariety of codimension $m+1$ in Gromov-Hausdorff sense. We also show that the K\"ahler-Ricci flow resolves certain type of conical singularities in Gromov-Hausdorff sense.

\section{Introduction}

The formation of singularities of the K\"ahler-Ricci flow on a compact manifold $M$ reveals the analytic and algebraic structures of $M$. It is well known that the K\"ahler-Ricci flow converges to a K\"ahler-Einstein metric if $M$ admits negative or vanishing first Chern class for any initial K\"ahler metric \cite{Cao1}.

When the canonical bundle $K_M$ is not nef, the K\"ahler-Ricci flow will develop finite time singularity and one expects the Ricci flow to carry out surgeries through the singularities in some natural and unique way.  The flow $(M, g(t))$ should converge in some suitable sense to a  `limit manifold'  $(\ov{M}, g_T)$ as $t$ tends to the singular time $T$  and continue on the new manifold starting at $g_T$. This is referred to as \emph{canonical surgery by the Ricci flow}. If the K\"ahler manifold $M$ is projective, then one  hopes that the canonical surgeries correspond to algebraic transformations such as divisorial contractions or flips.

An analytic analogue of Mori's minimal model program is laid out in [SoT3] for how the K\"ahler-Ricci flow will behave on a general projective variety.  More precisely, it is conjectured that the K\"ahler-Ricci flow will either deform a projective variety $M$ to its minimal model after finitely many divisorial contractions and flips in Gromov-Hausdorff sense, or collapse in finite time.  The existence and uniqueness is proved in \cite{SoT3} for the weak solution of the K\"ahler-Ricci flow through divisorial contractions and flips.  However, the Gromov-Hausdorff convergence at the singular time is largely open.  The program is established for K\"ahler surfaces in \cite{SW2, SW3}. More precisely, for the K\"ahler-Ricci flow on a K\"ahler surface $M$ with initial K\"ahler metric $g_0$,  either the flow deforms $M$ to a minimal surface  or the volume tends to $0$  in finite time,  after finitely many  contractions of $(-1)$-curves in Gromov-Hausdorff sense.   In the work of \cite{LT} the conjectural behavior of the flow through a flip is discussed in relation to their V-soliton equation.

The goal of the current paper is to establish examples of small contractions and resolution of singularities by the K\"ahler-Ricci flow.  Let $(X, g_0)$ be a compact K\"ahler manifold of complex dimension $n \ge 2$. We write $\omega_0= \frac{\sqrt{-1}}{2\pi} (g_0)_{i \ov{j}}dz^i \wedge d\ov{z^j}$ for the K\"ahler form associated to $g_0$.  We consider the following K\"ahler-Ricci flow $\omega=\omega(t)$ given by
\begin{equation} \label{krf0}
\frac{\partial}{\partial t} \omega = - \textrm{Ric}(\omega), \quad \omega|_{t=0} = \omega_0,
\end{equation}
for $ \textrm{Ric}(\omega) = - \ddbar \log \omega^n$, and $g=g(t)$ is the metric associated to $\omega$.  The flow admits a smooth solution if and only if, the K\"ahler class of $\omega(t)$ is given by
$$[\omega(t)] = [\omega_0] + t  [K_X] >0.$$ The first singular time  $T$ is characterized by
\begin{equation} \label{T}
T = \sup \{ t \in \mathbb{R} \ | \ [\omega_0] + t  [K_X] >0 \}.
\end{equation}
Clearly $T$ depends  only on $X$ and the K\"ahler class $[\omega_0]$, and satisfies  $0 < T \le \infty$.

The manifold $$X_{m,n}= \mathbb{P}(\mathcal{O}_{\mathbb{P}^n} \oplus \mathcal{O}_{\mathbb{P}^n}(-1)^{\oplus (m+1)}) $$ is a projective toric manifold for $m\geq 0$ and $n\geq 1$. $X_{0,n}$ is exactly $\mathbb{P}^{n+1}$ blown up at one point. $X_{m,n}$ does not admit a definite or vanishing first Chern class when $n\leq m$ and $X_{m,n}$ is Fano if and only if $n>m$. $X_{m,n}$ has a special subvariety $P_0$ of codimension $m+1$, defined as the zero section of projection $\mathbb{P}(\mathcal{O}_{\mathbb{P}^n} \oplus \mathcal{O}_{\mathbb{P}^n}(-1)^{\oplus (m+1)})\rightarrow \mathbb{P}^n$.  There exists a morphism
\begin{equation}
\Phi_{m,n}: X_{m,n} \rightarrow \mathbb{P}^{(m+1)(n+1)}
\end{equation} which is an immersion on $X_{m,n} \setminus P_0$ and contracts $P_0$ to a point. $Y_{m,n}$, the image of $X_{m,n}$ via $\Phi_{m,n}$ is smooth if and only if $m=0$ and then $Y_{0,n}$ is simply $\mathbb{P}^{n+1}$. When $m\geq 1$, $Y_{m,n}$ has a conical singularity at where $P_0$ is contracted. In particular, $Y_{m,n}=Y_{n,m}$ is the projective cone in $\mathbb{P}^{(m+1)(n+1)}$ over $\mathbb{P}^{m}\times \mathbb{P}^n$ via the Segre map for $m\geq 1$. It is also well-known that $X_{m,n}$ and $X_{n,m}$ are birationally equivalent for $m\geq 1$, and differ by a flip for $m\neq n$.

In this paper, we always consider the K\"ahler metrics on $X_{m,n}$ satisfying the Calabi symmetric condition defined in \cite{C1}. The precise definition is given in section \ref{calabi}.

Our first main result characterizes the limiting behavior of the K\"ahler-Ricci flow (\ref{krf0}) on $X_{m,n}$ as $t\rightarrow T$.

\begin{theorem} \label{main1} Let $g(t)$ be the  solution of the K\"ahler-Ricci flow (\ref{krf0}) on $X_{m,n}$ with initial K\"ahler metric $\omega_0 \in a_0 [D_H] + b_0 [D_\infty] $ satisfying the Calabi symmetry. Let $T>0$ be the first singular time of the flow.

\begin{enumerate}

\item If $m<n$ and $b_0/(m+2) > a_0/(n-m)$, then $T= a_0/(n-m)$ and  on $X_{m,n}\setminus P_0$,  $g(t)$ converges smoothly to a K\"ahler metric $g_T$.   Let $(X_T, d_T)$ be the metric completion of $( X_{m,n} \setminus P_0,g_T)$. Then $(X_T, d_T)$ has finite diameter and is homeomorphic to $Y_{m,n}$ as the projective cone in $\mathbb{P}^{(m+1)(n+1)}$ over $\mathbb{P}^{m}\times \mathbb{P}^n$ via the Segre map. Furthermore, $(X_{m,n}, g(t))$ converges to $(X_T, d_T)$ in Gromov-Hausdorff sense as $t\rightarrow T$.

\item If $m<n$ and $b_0/(m+2) = a_0/(n-m)$, then $T=a_0/(n-m)$ and $(X_{m,n}, g(t))$ converges to a point in Gromov-Hausdorff sense as $t\rightarrow T$.

\item If $m<n$ and $b_0/(m+2) < a_0/(n-m)$, then $T=b_0/(m+2)$ and $(X_{m,n}, g(t))$ converges to $(\mathbb{P}^n, (a_0-\frac{n-m}{m+2}b_0)\omega_{FS})$  in Gromov-Hausdorff sense as $t\rightarrow T$, where $\omega_{FS}$ is the Fubini-Study metric on $\mathbb{P}^n$.

\item If $m\geq n$,  then $T=b_0/(m+2)$ and $(X_{m,n}, g(t))$ converges to $(\mathbb{P}^n, (a_0-\frac{n-m}{m+2}b_0)\omega_{FS})$  in Gromov-Hausdorff sense as $t\rightarrow T$.

\end{enumerate}

\end{theorem}

In the case (3) and (4), the K\"ahler-Ricci flow can be continued on $\mathbb{P}^n$ starting with $(\mathbb{P}^n, (a_0-\frac{n-m}{m+2}b_0)\omega_{FS})$ and the flow will eventually become extinct in finite time. The case (2) is related to the result in \cite{So} that the K\"ahler-Ricci flow shrinks to a point if and only if $X$ is Fano and the initial K\"ahler class is proportional to $c_1(X)$, establishing the smooth case of a conjecture in \cite{T2}.  If this occurs, it is natural to renormalize the flow so that the volume is constant.  The problem of how  this normalized  flow behaves  is related to various notions of stability \cite{Y2,T2,Do} and is still open in general.  Assuming the existence of a K\"ahler-Einstein metric \cite{P2,TZhu} or soliton \cite{TZhu}, the flow is shown to converge to a K\"ahler-Einstein metric or soliton respectively (see also \cite{SeT, Zhu}).  The connection between stability conditions and the behavior of the K\"ahler-Ricci flow has been studied in \cite{PS, PSSW1, PSSW2, R, Sz, To, MS, CW} for example.

Theorem \ref{main1} can also be viewed as an analogue of Theorem 1.1 in \cite{SW1}. We apply ideas and techniques from \cite{SW1} to obtain many estimates in the proof of Theorem \ref{main1}. In fact, Theorem \ref{main1} can be generalized to  $$X_{m,n, k}= \mathbb{P}(\mathcal{O}_{\mathbb{P}^n} \oplus \mathcal{O}_{\mathbb{P}^n}(-k)^{\oplus (m+1)}),~~~k=1, 2, ... $$ In particular, $X_{m,n,1}= X_{m,n}$ and $X_{0,1,k}$ are exactly the rational ruled surfaces considered in \cite{SW1}.

We would also like to mention some known results about K\"ahler-Ricci solitons on these manifolds. $X_{n,0}$ admits a K\"ahler-Ricci soliton \cite{Koi, Cao2}. Complete K\"ahler-Ricci solitons are also constructed on vector bundles $ \mathcal{O}_{\mathbb{P}^n}(-1)^{\oplus (m+1)}$ by \cite{FIK} when $m=0$ and by \cite{Li} when $m\geq 1$.

The general conjecture in \cite{SoT3} predicts that the flow can also be continued in the first case in Theorem \ref{main1} and the contracted variety should jump to its minimal resolution by the flip. Our next result shows how the K\"ahler-Ricci flow can resolve certain type of projective conical singularities and confirms the weaker statement of the general conjecture.

\begin{theorem} \label{main2} Let $Y_{m,n}=Y_{n,m}$ be a projective cone in $\mathbb{P}^{(m+1)(n+1)}$ over $\mathbb{P}^m\times \mathbb{P}^n$ via the Segre map and let $g_0 \in [\mathcal{O}_{\mathbb{P}^{(m+1)(n+1)}}(1)]$ be the restriction of the Fubini-Study metric of $\mathbb{P}^{(m+1)(n+1)}$ on $Y_{m,n}$.

\begin{enumerate}

\item If $m>n\geq 1$, there exists a  smooth solution $g(t)\in \Phi_{m,n}^*[\mathcal{O}_{\mathbb{P}^{(m+1)(n+1)}}(1)]+t[K_{X_{m,n}}]$ of the K\"ahler-Ricci flow on $(0, T=1/(m+2))\times X_{m,n}$ such that on $X_{m,n}\setminus P_0\simeq Y_{m,n}\setminus\{O\}$, $g(t)$ converges smoothly to $g_0$ as $t\rightarrow 0$, and $(X_{m,n}, g(t))$ converges to $(Y_{m,n}, g_0)$ in  Gromov-Hausdorff sense as $t\rightarrow 0$. Furthermore, $g(t)$ has uniformly bounded local potential in $L^\infty$ for $t\in [0, T)$. If there exists another solution $\hat g(t)$ satisfying the above conditions, then $g(t)=\hat g(t)$.

\item If $m=n\geq 1$, there exists a  smooth solution $g(t) \in (1-(m+2)t) [g_0]$ of the K\"ahler-Ricci flow on $Y_{n,n}\setminus\{O\}$ for $t\in (0, T=1/(m+2) )$  such that

\begin{itemize}

\item $(Y_t, d_t)$, the metric completion of $(Y_{n,n}\setminus \{O\}, g(t))$ is homeomorphic to $(Y_{n,n}, g_0)$.

\item $g(t)$ converges to $g_0$ smoothly on $Y_{n,n}\setminus \{O\}$ as $t\rightarrow 0$, and $g(t)$ has uniformly bounded local potential in $L^\infty$ for $t\in [0, T)$.

\item  $(Y_{n,n}, d_t)$ converges to $(Y_{n,n}, g_0)$ in Gromov-Hausdorff sense as $t\rightarrow 0$ and converges to a point in Gromov-Haudorff sense as $t \rightarrow T$.

\end{itemize}

If there exists another solution $\hat g(t)$ satisfying the above conditions, then $g(t)=\hat g(t)$.

\end{enumerate}

\end{theorem}

The above theorem shows that the K\"ahler-Ricci flow resolves the conical singularity of $Y_{m,n}$ for $m>n\geq 1$ in Gromov-Hausdorff sense. It suggests that the K\"ahler-Ricci flow  smoothes out not only the initial singular metric, but also the initial underlying variety. Combining $(1)$ in Theorem \ref{main1} and $(1)$ in Theorem \ref{main2},  the K\"ahler-Ricci flow replaces $X_{m,n}$ by $X_{n,m}$ as an analytic flip in Gromov-Hausdorff sense if we are allowed to continue the K\"ahler-Ricci flow at $t=T$ with the Fubini-Study metric restricted on $Y_{m,n}$. We believe that the K\"ahler-Ricci flow should perform the flip for $X_{m,n}$ without replacing the singular metric $g_T$ by the Fubini-Study metric as the initial metric at the singular time.

On the other hand, the K\"ahler-Ricci flow does not change the underlying manifold $Y_{n,n}$ for $n\geq 1$ even though $K_{Y_{n,n}}$ is not a Cartier divisor. This is because $X_{n,n}$ is the resolution of $Y_{n,n}$ and the canonical divisor $K_{X_{n,n}}$ vanishes along the exceptional locus over the singularity $O$.

The organization of the paper is as follows. In section 2, we describe flips, the Calabi ansatz and the K\"ahler-Ricci flow on $\mathbb{P}(\mathcal{O}_{\mathbb{P}^n} \oplus \mathcal{O}_{\mathbb{P}^n}(-1)^{\oplus (m+1)})$. In section 3, we prove the small contraction by the K\"ahler-Ricci flow if the volume does not tend to zero when approaching the singular time. In section 4, we show that the K\"ahler-Ricci flow collapses if the volume tends to zero when approaching the singular time. In section 5, we describe how the K\"ahler-Ricci flow resolves singularities of $Y_{m,n}$.

\section{Background}

\subsection{An example of flips}

We will describe a family of projective bundles over $\mathbb{P}^n$ so that one can construct a flip. The detailed algebraic construction can be found in section 1.9 of \cite{D}.

Let $E$ be the vector bundle  over a projective space $\mathbb{P}^n$ defined by $ \mathcal{O} _{\mathbb{P}^n} \oplus \mathcal{O}_{\mathbb{P}^n} (-1)^{\oplus (m+1)}$.  We let
$$X_{m, n} = \mathbb{P} ( \mathcal{O}_{\mathbb{P}^n}  \oplus \mathcal{O}_{\mathbb{P}^n} (-1)^{\oplus (m+1)})$$
be the projectivization of $E$ and it is a $\mathbb{P}^{m+1}$ bundle over $\mathbb{P}^n$.  In particular, $X_{0, n}$ is $\mathbb{P}^{n+1}$ blown up at one point. Let $D_\infty$ be the divisor in $X_{m,n}$ given by $\mathbb{P}(\mathcal{O}_{\mathbb{P}^n}(-1)^{\oplus (m+1)})$, the quotient of $\mathcal{O}_{\mathbb{P}^n}(-1)^{\oplus (m+1)}$. We also let $D_0$ be the divisor in $X_{m,n}$ given by  $\mathbb{P}(\mathcal{O}_{\mathbb{P}^n} \oplus \mathcal{O}_{\mathbb{P}^n} (-1) ^{\oplus m} )$, the quotient of $\mathcal{O}_{\mathbb{P}^n} \oplus \mathcal{O}_{\mathbb{P}^n} (-1) ^{\oplus m} $. In fact, $N^1(X_{m,n})$ is spanned by
$[D_0]$ and $[D_\infty]$. We also define the divisor $D_H$ on $X_{m,n}$ by the pullback of the divisor on $\mathbb{P}^n$ associated  to $\mathcal{O}_{\mathbb{P}^n}(1)$.
Then $$[D_\infty]= [D_0] + [D_H]$$ and
\begin{equation}\label{canbundle}
[ K_{X_{m,n}} ] = - (m+2) [D_\infty] - (n-m) [D_H] = -(n+2) [D_\infty] + (n-m) [D_0].
\end{equation}
The above formulas can be easily obtained by induction on $m$ and the adjunction formula.
In particular, $D_\infty$ is a big and semi-ample divisor and any divisor $ a[D_H] + b[D_\infty]$ is ample if and only if $a>0$ and $b>0$. Hence $X_{m,n}$ is Fano if and only if $n>m$.

Let  $P_0$ be  the zero section of $\pi_{m.n}: X_{m,n} \rightarrow \mathbb{P}^n$, which is the intersection of the $m+1$ effective divisors as the quotient of $\mathcal{O}_{\mathbb{P}^n} \oplus \mathcal{O}_{\mathbb{P}^n} (-1) ^{\oplus m}$. In fact, the linear system $| [D_\infty]|$ is base-point-free and it induces a morphism

$$\Phi_{m,n}: X_{m,n} \rightarrow \mathbb{P}^{(m+1)(n+1)}.$$ $\Phi_{m,n}$ is an immersion on $X_{m,n}\setminus P_0$ and it contracts $P_0$ to a point. $Y_{m,n}$, the image of $\Phi_{m,n}$ in $\mathbb{P}^{(m+1)(n+1)}$, is a projective cone over $\mathbb{P}^m\times \mathbb{P}^n $ in $\mathbb{P}^{(m+1)(n+1)}$ by the Segre embedding $$[Z_0, ..., Z_m]\times[W_0, ..., W_n]\rightarrow [Z_0W_0, ..., Z_iW_j, ..., Z_mW_n]\in \mathbb{P}^{(m+1)(n+1)-1} .$$

\subsection{Calabi ansatz}\label{calabi}

In this section, we will define the Calabi ansatz constructed by Calabi \cite{C1} (also see \cite{Li}). To apply the Calabi symmetry, we instead consider the vector bundle $$E= \mathcal{O}_{\mathbb{P}^n}(-1) ^{\oplus (m+1)}.$$

Let $\omega_{FS}$ be the Fubini-Study metric on $\mathbb{P}^n$. Let $h$ be the hermitian metric on $\mathcal{O}_{\mathbb{P}^d} (-1)$ such that $Ric(h) = -\omega_{FS}$. The induced hermtian metric $h_E$ on $E$ is given by $h_E= h^{\oplus (m+1)}$.
Under local trivialization of $E$, we write
$$e^\rho = h_\xi (z) |\xi|^2, ~~~ \xi = (\xi_1, \xi_2, ..., \xi_{m+1}),$$ where $h_\xi(z)$ is a local representation for $h$ (note that $h_E$ has same eigenvalues $h(z)$). In particular, if we choose the inhomogeneous coordinates $z=(z_1, z_2, ..., z_n)$ on $\mathbb{P}^n$, we have
$$h_\xi (z)= (1+|z|^2).$$
We would like to find appropriate conditions for $a\in \mathbb{R}$ and $u(\rho)$ such that
\begin{equation}\label{metricrep1}
\omega = a \omega_{FS} + \ddbar u(\rho)
\end{equation}
defines a K\"ahler metric on $X_{m,n}$.
In fact,
\begin{equation}\label{metricrep2}
\omega = (a + u'(\rho)) \omega_{FS} + \frac{\sqrt{-1}}{2\pi} h_\xi e^{-\rho} ( u' \delta_{\alpha \beta} + h_\xi e^{-\rho} ( u'' - u') \xi^{\bar \alpha} \xi^{\beta} ) \nabla \xi^\alpha \wedge \nabla \xi^{\bar \beta}.
\end{equation}
Here, $$\nabla \xi^\alpha = d \xi^\alpha + h_\xi^{-1} \partial h_\xi \xi^\alpha$$ and $\{ dz^i, \nabla \xi^\alpha\}$ is dual to the basis

$$ \nabla_{z^i} = \frac{\partial}{\partial z^i} - h_\xi ^{-1} \frac{\partial h_\xi }{\partial z^i} \sum_\alpha \xi^\alpha \frac{\partial }{\partial \xi^\alpha}, ~~ \frac{\partial}{\partial \xi ^\alpha}.$$

The following criterion is due to Calabi \cite{C1}.
\begin{proposition}\label{kacon}

$\omega$ as defined above is a K\"ahler metric if and only if

\begin{enumerate}

\item

$a>0$.

\item $u'>0$ and $u''>0$ for $\rho\in (-\infty, \infty)$.

\item $U_0 (e^\rho) = u(\rho)$ is smooth on $(-\infty, 0]$ and $U_0' (0)>0$.

\item $U_\infty (e^{-\rho}) = u(\rho) - b \rho$ is smooth on $[0, \infty)$ for some $b>0$ and $U_\infty'(0)>0$.

\end{enumerate}

\end{proposition}

We remark that given $a$, $b>0$, the K\"ahler metric constructed above lies in the K\"ahler class

\begin{equation}
\omega=a\omega_{FS} + \ddbar u (\rho) \in a[D_H]+ b[D_\infty]
\end{equation}
and \begin{equation}\label{autobd} 0< u'(\rho)\leq b.\end{equation}

\subsection{The K\"ahler-Ricci flow on $X_{m,n}$}\label{redric}

Straightforward calculations show that the induced volume form of $\omega$ is given by

\begin{equation}\label{volrep}
\omega^{m+n+1} = (a + u')^n h_\xi^{m+1} e^{-(m+1)\rho} ( u')^{m} u''  (\omega_{FS}^n \wedge \prod_{\alpha=1}^{m+1} \frac{\sqrt{-1}}{2\pi} d\xi^\alpha\wedge d \xi^{\bar \alpha}).
\end{equation}
Therefore,

$$- Ric(\omega) = \ddbar (  \log [(a+ u')^n (u')^{m} u''] - (m+1)\rho )  +  ( m  - n  )\omega_{FS}.$$

It is straightforward to check that  the Calabi ansatz is preserved by the Ricci flow. Indeed, the K\"ahler-Ricci flow
\begin{eqnarray}
 \ddt{\omega} = - Ric (\omega),  ~~~\omega|_{t=0} = \omega_0= a_0 \omega_{FS} + \ddbar u_0
\in a_0[D_H]+b_0[D_\infty]
\end{eqnarray} is equivalent to the following parabolic equation
$$ a'(t)\omega_{FS} + \ddbar \ddt{u} = (m-n) \omega_{FS} + \ddbar (  \log [(a+ u')^n (u')^{m} u''] - (m+1) \rho ) .
 $$
Separating the variables, we have that

\begin{equation}
a=a(t) =a_0  - ( n-m) t
\end{equation}
 and
\begin{equation}\label{krfu}
\ddt{u} = \log [(a+ u')^n (u')^{m} u''] - (m+1)\rho + c_t,
\end{equation}
where
\begin{equation} \label{ct}
c_t = - \log u''(0,t) - m \log u'(0,t)- n\log (a(t)+u'(0,t)).
\end{equation}
From the formula (\ref{canbundle}) and the K\"ahler class evolves by $[\omega]=(a_0-(n-m)t)[D_H]+ (b_0-(m+2)t) [D_\infty]$, and so $$ b=b(t) = b_0 - (m+2) t. $$

 It is straightforward to show that equation (\ref{krfu}) admits a smooth solution $u$ satisfying the Calabi ansatz as long as the K\"ahler-Ricci flow admits a smooth solution, by comparing $u$ to the solution of the Monge-Amp\`ere flow associated to the K\"ahler-Ricci flow.

Next,  the evolution equations for $u'$ and $u''$ are given by
\begin{eqnarray}
 \ddt{u'} &=& \frac{u'''}{u''} + \frac{ m u''}{u'} + \frac{ n u''}{a + u'} - (m+1),
 \\ \label{udpevolution}
\ddt{u''}& =& \frac{u^{(4)}} {u''} - \frac{ (u''')^2 }{ (u'')^2} + \frac{ m u'''}{u'} - \frac{ m (u'')^2}{ (u')^2} + \frac{ n u'''}{ a + u'} - \frac{ n (u'')^2}{( a + u')^2},
\end{eqnarray}
as can be seen from differentiating (\ref{krfu}).

\section{Small contractions by the K\"ahler-Ricci flow}

The first singular time of the K\"ahler-Ricci flow on $X_{m,n}$ is given by
\begin{equation}
T=\sup\{ t>0~|~[\omega_0] + t [K_{X_{m,n}}] >0\}.
\end{equation}
Since $X_{m,n}$ is not a minimal model, $T<\infty$.

In the section, we assume that at the singular time $T$, $a(T)=0$ and $b(T)>0$, i.e., the K\"ahler-Ricci flow does not collapse. This is equivalent to,

$$ n>m, ~~\frac{b_0}{m+2} > \frac{a_0}{n-m}.$$
In this case, the first singular time of the flow is given by $T = \frac{a_0}{n-m}.$

Let us first explicitly write down the contraction map $\Phi_{m,n}$. In local trivialization for $X_{m,n}$,  $\{ 1, \xi_{\alpha}, z_i \xi_{\alpha} \}_{ i=1, ..., n, \alpha = 1, ..., m+1}  $ extend to global holomorphic sections in $[D_\infty]$, furthermore, they span
$H^0(X_{m,n}, \mathcal{O}([D_\infty])) $. Then the free linear system of $|[D_\infty]|$ induces the following morphism
$$ \Phi_{m,n}: (z_i, \xi_\alpha) \in X_{m,n} \rightarrow [1, \xi_\alpha, z_i \xi_{\alpha} ] \in \mathbb{P}^{(m+1)(n+1)}.$$
The pullback of the Fubini-Study metric is given by
\begin{eqnarray*}\hat\omega &=& \ddbar \log ( 1+ \sum_{\alpha=1}^{m+1}|\xi_\alpha|^2 + \sum_{1 \leq i \leq n, 1 \leq \alpha \leq m+1} |z_i \xi_\alpha|^2 )\\
&=& \ddbar\log (1+ (1+\sum_{i=1}^n |z_i|^2) (\sum_{\alpha=1}^{m+1} |\xi_\alpha|^2)) \\
&=&\ddbar\log (1+ e^\rho) .
\end{eqnarray*}

 Let   $$\hat u (\rho) = \log (1 + e^\rho). $$
Then $\ddbar \hat u$ extends to the pullback of the Fubini-Study metric $\hat\omega$ given by $\Phi_{m,n}$. In particular, $Y_{m,n}$ has an isolated conical singularity and $\hat\omega$ is a asymptotically conical metric on $Y_{m,n}$ near the conical singularity.

Now we list some well-known results for some useful uniform estimates. We begin by rewriting the K\"ahler-Ricci flow as a parabolic flow of Monge-Ampere type. We let $\omega_0$ be the initial K\"ahler metric and $\Omega$ a smooth volume form on $X_{m,n}$. Let $\chi = \ddbar \log \Omega$ and $\omega_t = \omega_0 + t\chi \in [\omega_0]+t[K_{X_{m,n}}]$ be the reference form. Then the K\"ahler-Ricci flow is equivalent to

\begin{equation}\label{0ord1}
\ddt{\varphi} = \log \frac{ (\omega_t + \ddbar \varphi)^{m+n+1}}{\Omega} , ~~~~\varphi|_{t=0}=0.
\end{equation}

Then there exists a unique solution $\varphi\in C^\infty([0, T)\times X_{m,n} )$. Furthermore,  we have the following well-known estimates due to \cite{TZha}.

\begin{enumerate}

\item There exists $C>0$ such that on $[0, T)\times X_{m,n}$, \begin{equation}| \varphi |  \leq C. \end{equation}

\item There exists $C>0$ such that on $[0, T)\times X_{m,n}$, \begin{equation}\label{volest} (\omega_t + \ddbar\varphi)^{m+n+1} \leq C \Omega.\end{equation}

\item  For any $K\subset\subset X_{m,n}\setminus P_0$, there exists $C>0$ such that on $[0, T)\times X_{m,n}$, \begin{equation}\label{higherord1} || \varphi ||_{C^k ([0, T) \times K)} \leq C_{k, K}. \end{equation}

\end{enumerate}

We also have the following estimates as the parabolic Schwarz lemma.
\begin{lemma} There exists $C>0$ such that on $[0, T)\times X_{m,n}$,

$$\omega \geq C \hat \omega.$$

\end{lemma}

\begin{proof} The proof is given in \cite{SW1, So} and makes  use of the $L^\infty$-estimate of $\varphi$.

\qed
\end{proof}

By comparing $\omega(t)$ and $\hat \omega$, we immediately have the following estimate.

\begin{corollary} There exists $C>0$ such that on $[0, T)\times X_{m,n}$,

$$ a+ u' \geq C \hat u' = \frac{C e^\rho}{(1+e^\rho)} .$$

\end{corollary}

\begin{lemma} There exists $C>0$ such that on $[0, T)\times X_{m,n}$,

\begin{equation} \label{vol3}((u')^{m+n+1})' \leq \frac{C e^{(m+1)\rho}}{(1+e^\rho)^{m+2}} .\end{equation}

\end{lemma}

\begin{proof}  The inequality (\ref{vol3}) follows immediately from the volume estimate (\ref{volest}) by the following observation
$$((u')^{m+n+1})' \leq (a+ u')^n (u')^{m} u'' \leq C_1 (1+ \hat u')^n (\hat u')^{m} \hat u''\leq  \frac{C_2e^{(m+1)\rho}}{(1+e^\rho)^{m+2}}. $$
The last inequality follows from the definition of $\hat u$.

\qed
\end{proof}

\begin{corollary}\label{1stder1}  There exists $C>0$ such that on $ [0, T) \times X_{m,n}$,

\begin{equation}
 u'(\rho) \leq C e^{\frac{m+1}{m+n+1} \rho}.
\end{equation}

\end{corollary}

\begin{proof} Using Proposition \ref{kacon} and integrating (\ref{vol3}) from $-\infty$ to $\rho$, we have

$$(u'(\rho))^{m+n+1} \leq C  \int_{-\infty}^\rho  e^{(m+1)\rho} d\rho + \lim_{\rho\rightarrow -\infty} (u'(\rho))^{m+n+1} =  \frac{C}{m+1} e ^{(m+1)\rho}.$$

\qed
\end{proof}

We also notice that $0< u'(\rho) < b(t)$ for $\rho\in(-\infty, \infty)$ because $u'$ is increasing and $\lim_{\rho\rightarrow \infty} u'(\rho)=  b(t)$.  Therefore, $u'$ is uniformly bounded above for $t\in [0, T)$.

\begin{proposition}\label{2ndder1}There exist $C>0$ such that on $ [0, T)\times X_{m,n}$,

\begin{equation}
u'' \leq C u' .
\end{equation}

\end{proposition}

\begin{proof}

Let $H= \log u'' - \log u'$. Notice that by Proposition \ref{kacon}, for fixed $t\in[0, T)$ and near $\rho=-\infty$, $u(\rho) = U_0(e^{\rho})$ for some smooth function $U_0$, and near $\rho=\infty$,  $u(\rho) = U_\infty(e^{-\rho})+ b \rho$ for some smooth function $U_\infty$ and $b>0$.    $$\
\lim_{\rho\rightarrow -\infty} \frac{u''}{u'}  = \lim_{\rho\rightarrow -\infty}  \frac{ (U_0(e^\rho))''}{ (U_0(e^\rho))'}=\lim_{\rho\rightarrow -\infty}  (  \frac{U_0'  +  e^{\rho}  U_0''} {  U_0' }) = 1 + \lim_{\rho\rightarrow -\infty} e^\rho \frac{U_0''}{U_0'} =1.$$

$$\lim_{\rho\rightarrow \infty}  \frac{u''}{u'}  =  \lim_{\rho\rightarrow \infty} \frac{ (U_\infty(e^{-\rho}) + b\rho)''}{ (U_\infty(e^{-\rho}) + b\rho)'}= \lim_{\rho\rightarrow \infty} \frac{e^{-\rho} U_\infty'  +  e^{-2\rho}  U_\infty'' } { -e ^{-\rho} U_\infty' + b} =0.$$ And so we can apply maximum principle for $H$ in $[0, T)\times (-\infty, \infty)$.

\begin{eqnarray*}
\ddt{H}
&=&  \frac{1}{u''} \{ \frac{ u^{(4)} }{u''} - \frac{ (u''')^2}{ (u'')^2} + \frac{  mu'''}{ u'} - \frac{ m (u'')^2}{ (u')^2} + \frac{ n u'''}{a+ u'} - \frac{ n (u'')^2}{ (a + u')^2} \} \\
&& - \frac{1}{u'} \{ \frac{u'''}{u''} + \frac{ m u''}{u'} + \frac{ n u''}{ a+ u'} - (m+1)\}
\end{eqnarray*}

Suppose that $H(t_0, \rho_0) = \sup_{[0, t_0] \times (-\infty, \infty) } H(t, \rho)$ is achieved for some $t_0\in (0, T)$, $\rho_0\in (-\infty, \infty)$.  At $(t_0, \rho_0)$, we have

$$H'= \frac{ u'''}{u''} - \frac{ u''}{u'} = 0$$ and

$$ H''= \frac{ u^{(4)}} {u''} - \frac{ (u''')^2}{(u'')^2} - \frac{ u''' } { u'} + \frac{ (u'')^2}{ (u')^2} = \frac{ u^{(4)} } {u''} - \frac{ u''' } { u'}\leq 0. $$
Then at $(t_0,\rho_0)$,
\begin{eqnarray*}
0 &\leq&\ddt{H}\\
&= & \frac{1}{u''} \{ \frac{ u^{(4)} }{u''} - \frac{ (u''')^2}{ (u'')^2} + \frac{ m u'''}{ u'} - \frac{ m (u'')^2}{ (u')^2}  \}  - \frac{1}{u'} \{ \frac{u'''}{u''} + \frac{ m u''}{u'}  - (m+1)\}  \\
&& + \frac{n}{a+u'} \{  \frac{u'''}{u''} - \frac{u''}{a+u'} - \frac{u''}{u'} \}  \\
&\leq&  - \frac{(m+1) u''}{(u')^2} + \frac{m+1}{u'}  - \frac{n u''}{(a+u')^2}\\
&\leq&  \frac{m+1}{u'} (1- e^H ).
\end{eqnarray*}
Therefore by the maximum principle, $H(t_0, x_0) \leq 0$ and so
$$\sup_{[0, T)\times (-\infty, \infty)} H(t, \rho) \leq \sup_{(-\infty, \infty)} H(0, \rho) < \infty.$$
The proposition then follows .

\qed
\end{proof}

We have the following immediate corollary by combining Proposition \ref{2ndder1} and Corollary \ref{1stder1}.
\begin{corollary} \label{2ndder2}

There exists $C>0$ such that on $ [0, T) \times X_{m,n}$,

$$ u''(\rho) \leq C e^{\frac{m+1}{m+n+1} \rho}.$$

\end{corollary}

\begin{corollary}\label{metest} There exists $C>0$ such that  on $ [0, T) \times X_{m,n}$,

\begin{equation}
\omega \leq C (\omega_{FS} + \hat\omega+ e^{-\frac{n}{m+n+1}\rho} \hat\omega ).
\end{equation}

\end{corollary}

\begin{proof} The corollary holds for $t\in[0, T)$ and $\rho\in [0, \infty)$ from the estimates in (\ref{higherord1}) away from $P_0$ since $\omega_{FS}+\hat\omega$ is a smooth K\"ahler metric on $X_{m,n}$. It suffices to prove the corollary for $\rho\leq 0$.

Applying Corollary \ref{1stder1} and Corollary \ref{2ndder2}, for $(t, \rho)\in [0, T)\times(-\infty, 0]$, we have
\begin{eqnarray*}
\omega %
&\leq& C_1(1 + e^{\frac{m+1}{m+n+1} \rho}) \omega_{FS} + C_1 h_\xi e^{-\frac{n}{m+n+1} \rho} ( \delta_{\alpha \beta} + h_\xi e^{ -\rho} \xi^{\bar \alpha} \xi^{\beta} ) \frac{\sqrt{-1}}{2\pi} \nabla \xi^{\alpha} \wedge \nabla \xi^{\bar \beta} \\
& \leq& C_2 \omega_{FS} + C_2 h_\xi e^{-\frac{n}{m+n+1} \rho}  \sum_{\alpha=1}^{m+1}\frac{\sqrt{-1}}{2\pi}   \nabla \xi^{\alpha} \wedge \nabla \xi^{\bar \alpha}  .
\end{eqnarray*}
The corollary follows by comparing the above estimates to
$$\hat\omega= \frac{e^\rho}{1+e^\rho} \omega_{FS} + \frac{h_\xi}{1+e^\rho}  ( \delta_{\alpha \beta} - \frac{ h_\xi} {1+e^\rho} \xi^{\bar \alpha} \xi^{\beta} )  \frac{\sqrt{-1}}{2\pi}  \nabla \xi^{\alpha} \wedge \nabla \xi^{\bar \beta} \geq C_3 e^\rho \omega_{FS} + C_3 h_\xi \frac{\sqrt{-1}}{2\pi}   \nabla \xi^{\alpha} \wedge \nabla \xi^{\bar \alpha}  $$ for $\rho\leq 0$ and some $C_3>0$.

\qed
\end{proof}

Let $\omega(T)= \lim_{t\rightarrow T^-} \omega(t)$ be the closed positive $(1,1)$-form with bounded local potentials. Then by the estimates of $\omega(t)$ away from $P_0$ as in  (\ref{higherord1}), $\omega(T)$ is a smooth K\"ahler metric on $X_{m,n} \setminus P_0$ and $\omega(t)$ converges in $C^\infty(X_{m,n}\setminus P_0)$ to $\omega(T)$ as $t\rightarrow T$.

\begin{theorem}\label{ghcon2} Let $(X_T, d_T)$ be the metric completion of $( X_{m,n} \setminus P_0, \omega(T))$. Then $(X_T, d_T)$ has finite diameter and is homeomorphic to $Y_{m,n}$ as the projective cone in $\mathbb{P}^{(m+1)(n+1)}$ over $\mathbb{P}^{m}\times \mathbb{P}^n$ via the Segre map. Furthermore, $(X_{m,n}, g(t))$ converges to $(X_T, d_T)$ in Gromov-Hausdorff sense as $t\rightarrow T^-$, and there exists $C>0$ such that for $t\in [0, T)$.
$$diam( X_{m,n}, g(t)) \leq C.$$

\end{theorem}

\begin{proof} Let $U_\kappa = \{ e^\rho \leq \kappa\} $ be a $\kappa$-tubular neighborhood of the zero section $P_0$. We will use local coordinates $(z_i, \xi_\alpha)$ for $i=1, ..., n$ and $\alpha=1, ..., m+1$. For any fixed fibre $X_z= \pi^{-1}(z)$ for $z=(z_1, ..., z_n)\in \mathbb{C}^n$, there exists $C_1>0$, such that the restriction of the evolving metric is bounded by
\begin{eqnarray*}
\omega|_{X_z} &=& \frac{\sqrt{-1}}{2\pi} u''e^{-2\rho} \partial e^{\rho} \wedge \bar{\partial} e^{\rho} + u'e^{-\rho} \ddbar e^\rho - \frac{\sqrt{-1}}{2\pi} u' e^{-2\rho} \partial e^{\rho} \wedge \bar{\partial} e^{\rho}\\
&\leq& \frac{\sqrt{-1}}{2\pi} u''e^{-2\rho} \partial e^{\rho} \wedge \bar{\partial} e^{\rho} +  u'e^{-\rho}\ddbar e^\rho\\
%
&\leq& C_1 h_\xi e^{-\frac{n}{m+n+1}\rho} \sum_{\alpha} \frac{\sqrt{-1}}{2\pi}  d\xi^{\alpha} \wedge d\xi^{\bar \alpha}.
\end{eqnarray*}

We first show that for any $\epsilon>0$, there exists $\kappa_\epsilon>0$  such that for any $z\in \mathbb{C}^n$, $\kappa< \kappa_\epsilon$ and $t\in [0,  T)$,

$$diam(X_z\cap U_\kappa, g(t)) < \epsilon.$$

\begin{itemize}

\item We begin with estimates in the radial direction. We can always assume $\rho\leq 0$.
%
%
For any point $\xi\in X_z$, we consider the radial line segment $\gamma(r_\xi)$ joining $0$ and $\xi$ in $\mathbb{C}^{m+1}$ where $\xi= r_\xi e^{i \theta_\xi}$. Note that $e^\rho= (1+|z|^2)|\xi|^2$, then the arc length of $\gamma$ is given by
\begin{eqnarray*}
|\gamma|_{g(t)} &\leq& C_2 \int_0^{|\xi|} e^{-\frac{n}{2(m+n+1)}\rho} (1+|z|^2)^{1/2} d r\\
&=& C_2 \int_{0}^{|\xi|} (1+|z|^2) ^{\frac{m+1}{2(m+n+1)} } r^{ - \frac{n}{m+n+1}} d r \\
&\leq& C_3 \{ (1+|z|^2) r_\xi^2 \}^{\frac{m+1}{2(m+n+1)} }\\
&\leq& C_3 ~\kappa^{\frac{m+1}{2(m+n+1)}}  \end{eqnarray*}
for some fixed constant $C_i>0$, $i=2,3$.

\item We now consider the behavior of $g(t)$ on $S_{|\xi|}$, the sphere centered at $\xi=0$ with radius $r_\xi= |\xi|$ in $\mathbb{C}^{m+1}$ with respect to the Euclidean metric. Let $g_{S_{2m+1}}$ be the standard metric on the unit sphere $S_{2m+1}$ in $\mathbb{C}^{m+1}$. If $S_{|\xi|} \subset X_z\cap U_\kappa$, then there exist $C_4>0$ such that
\begin{eqnarray*}
g(t)|_{S_{|\xi|}} &\leq&  \sqrt{-1} C_1 e^{-\frac{n}{m+n+1}\rho} (1+|z|^2) d\xi \wedge d\bar\xi |_{S_{|\xi|}} \\
&\leq & C_4 e^{-\frac{n}{m+n+1}\rho} (1+|z|^2)  |\xi|^2 g_{S_{2m+1}}\\
& =& C_4 e^{\frac{m+1}{m+n+1} \rho} g_{S_{2m+1}}\\
&=& C_4 \kappa^{\frac{m+1}{m+n+1}} g_{S_{2m+1}}.
\end{eqnarray*}

\end{itemize}
Combining the above estimates, any two points in $X_z\cap U_\kappa$ can be connected with by a piecewise smooth curve in $X_z\cap U_\kappa$  with arbitrarily small arc length if $\kappa$ is chosen sufficiently small. 

Now we consider points $(0, \xi)$ and $(w, \xi) \in U_\kappa$ and we can assume  $\xi=(\xi_1, 0,..., 0)$ after the unitary transformation.  We will then consider a straight line segment $$\gamma(s)= \{ (z, \xi)~|~z=sw, ~\xi=(\xi_1, 0, ..., 0)\}.$$
There exists $C_5>0$, such that the restriction of $g(t)$ on the submanifold $V=\{ \xi= (\xi_1, 0, 0, ..., 0)\}$ is bounded by
\begin{eqnarray*}
g(t)|_{V\cap U_\kappa} &=& (a + u')\omega_{SF} + \frac{\sqrt{-1}}{2\pi} u''\frac{\bar z_i z_j}{(1+|z|^2)^2} \sum_{i,j} d z_i \wedge d z_{\bar j}\\
& \leq & C_5 (a + e^{\frac{m+1}{m+n+1}\rho})\omega_{SF} + \sqrt{-1} C_5 e^{\frac{m+1}{m+n+1}\rho} \frac{\bar z_i z_j}{(1+|z|^2)^2} \sum_{i,j} d z_i \wedge d z_{\bar j}\\
& \leq & C_5(a_0 -(n-m)t + \kappa^{\frac{m+1}{m+n+1}}) \omega_{SF} + \sqrt{-1} C_5 e^{\frac{m+1}{m+n+1}\rho} \frac{\bar z_i z_j}{(1+|z|^2)^2} \sum_{i,j} d z_i \wedge d z_{\bar j}.
\end{eqnarray*}

Therefore there exist $C_6, C_7>0$ such that the arc length of $\gamma(s)$ for $0\leq s\leq 1$ is bounded by
\begin{eqnarray*}
|\gamma| _{g(t)} &\leq& C_6 \int^1_0 e^{\frac{m+1}{2(m+n+1)}\rho} \frac{s|w|^2}{1+|z|^2} d s + C_6(a_0 -(n-m)t + \kappa^{\frac{m+1}{m+n+1}})^\frac{1}{2}\\
&\leq& C_6 \int^1_0 (1+|z|^2)^{\frac{m+1}{2(m+n+1)}-1}|\xi|^\frac{m+1}{m+n+1} s|w|^2 d s + C_6(a_0 -(n-m)t + \kappa^{\frac{m+1}{m+n+1}})^\frac{1}{2}\\
&\leq& C_7 \{(1+|w|^2)|\xi|^2\}^{\frac{m+1}{2(m+n+1)}} + C_7(a_0 -(n-m)t + \kappa^{\frac{m+1}{m+n+1}})^\frac{1}{2}\\
&\leq& C_7 \kappa^{\frac{m+1}{2(m+n+1)}} + C_7(a_0 -(n-m)t + \kappa^{\frac{m+1}{m+n+1}})^\frac{1}{2}.
\end{eqnarray*}

In general, given  two points $(z, \xi)$ and $(z', \xi') \in U_\kappa$, we can assume $|\xi|\leq |\xi'|$ without loss of generality. Let $\hat\xi=(\xi_1, 0,..., 0) $ such that $ |\hat\xi|=|\xi|$.
\begin{eqnarray*}
&&dist_{g(t)} ((z,\xi), (z',\xi')) \\
&\leq& dist_{g(t)} ((z,\xi), (z,\hat\xi)) +dist_{g(t)} ((z',\xi'), (z',\hat\xi))+ dist_{g(t)} ((0,\hat\xi), (z',\hat\xi)) +  dist_{g(t)} ((z,\hat\xi), (0,\hat\xi)) \\
&\leq& C_8 \kappa ^{\frac{n+1}{2(m+n+1)} }+ C_8 (a_0 - (n-m)t)^\frac{1}{2}
\end{eqnarray*}
for $C_8>0$. Hence for any $\epsilon>0$,  there exist $\kappa_\epsilon>0$ and $T_\epsilon \in (0, T)$ such that for any $0< \kappa < \kappa_\epsilon$ and $t\in (T_\epsilon, T)$,

$$ diam(U_\kappa\setminus P_0, g(t) ) < \epsilon.$$

This shows that $diam(X_{m,n}, g(t))$ is uniformly bounded above for $t\in[0, T)$. Similar argument shows that the metric completion $(X_T, d_T)$ of $(X_{m,n}\setminus P_0, \omega(T))$ is compact and is homeomorphic to $(Y_{m,n}, \hat g)$ as a metric space after replacing $a=a_0-(m-n)t$ by $0$.  Standard argument shows that $(X_{m,n}, g(t))$ converges to $(X_T, d_T)$ in  Gromov-Hausdorff sense as $t\rightarrow T$ (cf \cite{SW1}).

\qed
\end{proof}


\section{Finite time collapsing}\label{col}

\subsection{The case $m\geq n$,  or $m < n$ and $\frac{b_0}{m+2} <\frac{a_0}{n-m}$}\label{collapse}

In this section, we consider the K\"ahler-Ricci flow on $X_{m,n}$ with the initial class $a_0[D_H]+ b_0 [D_\infty]$ such that

$$m\geq n$$ or $$n>m,~~\frac{b_0}{m+2} < \frac{a_0}{n-m}.$$
The first singular time of the flow is given by \begin{equation} T=\frac{b_0}{m+2}.\end{equation} The following lemma is an immediate consequence of the observation (\ref{autobd}).

\begin{lemma} \label{obuprime}For $t\in [0, T)$ and $\rho\in(-\infty, \infty)$, we have

\begin{equation}
0< u' < b= b_0 - (m+2)t.
\end{equation}

\end{lemma}

The general volume estimate (\ref{volest}) gives us the upper bound for the volume form.

\begin{lemma}\label{volest2} There exists $C>0$ such that on $[0, T)\times X_{m,n}$,
$$\omega(t)^{m+n+1} \leq C \Omega.$$
\end{lemma}

\begin{corollary}
There exists $C>0$ such that for $\rho \in (-\infty, \infty)$ and $t\in [0, T)$,
\begin{equation}\label{haha1}
0<u' \leq C \min ( T-t, e^{ \rho})
\end{equation}
and
\begin{equation}\label{haha2}
0<b - u' \leq Ce^{-\frac{1}{m+1}\rho}.
\end{equation}

\end{corollary}

 \begin{proof} We apply similar argument as in Corollary \ref{vol3}.

\begin{itemize}

\item By Lemma \ref{volest2},
$$[ (u')^{m+1}]' \leq C_1(a+u')^n (u')^m u'' \leq C_2 (1+ \hat u')^n(\hat u')^m \hat u''. $$
For $\rho \in (-\infty, \infty)$,
$$  (u')^{m+1} (\rho) \leq C_3\int_{-\infty} ^\rho e^{(m+1)\rho} d \rho \leq C_4 e^{(m+1)\rho}.$$
The estimate (\ref{haha1}) follows by combining Lemma \ref{obuprime} and the same argument in Corollary \ref{1stder1}.

\item We also have $$[ (u')^{m+1}]' \leq C_5(a+u')^n (u')^m u'' \leq C_6 (1+ \hat u')^n(\hat u')^m \hat u''\leq C_7 e^{-\rho}.$$
Then
$$  b^{m+1} - (u')^{m+1} (\rho) \leq C_8\int^{\infty}_\rho e^{-\rho} d \rho \leq C_9 e^{-\rho}.$$
The estimate (\ref{haha1}) follows immediately.

\end{itemize}

\qed
\end{proof}

\begin{proposition} There exists $C>0$ such that on $[0, T)\times (-\infty, \infty)$,

\begin{equation}
u''\leq
C \min \{ u', b_t - u'\} .
\end{equation}

\end{proposition}

\begin{proof} The same argument in Proposition \ref{2ndder1} can be applied to show that $u''/u'$ is uniformly bounded above on $[0, T)\times (-\infty, \infty)$.

Let $H=\log \{u''/(b-u')\}$. Then $\lim_{\rho\rightarrow -\infty} H =-\infty $ and $\lim_{\rho\rightarrow \infty} H = 0$. The evolution equation for $H$ is given as follows.

\begin{eqnarray*}
\ddt{H} &=& \frac{1}{u''} \{ \frac{ u^{(4)}}{u''} - \frac{(u''')^2}{(u'')^2} + \frac{m u'''}{u'} - \frac{m(u'')^2}{(u')2} + \frac{n u'''}{a+u'} - \frac{n (u'')^2}{(a+u')^2}    \} \\
&&+ \frac{1}{b-u'} \{ \frac{u'''}{u''} + \frac{mu''}{u'} + \frac{n u''}{a+u'} - (m+1) \}.
\end{eqnarray*}

We also have
$$ H'=\frac{u'''}{u''} + \frac{u''}{b- u'}  $$ and
$$ H''=\frac{u^{(4)}}{u''} - \frac{ (u''')^2}{(u'')^2} + \frac{u''}{b-u'} H'   .$$

Suppose $\sup_{[0, t_0)\times (-\infty,\infty)} H = H(t_0, \rho)$. Then at $(t_0, \rho_0)$, straightforward calculations show that
$$0\leq \ddt{H} <0,$$ which is a contradiction. Thus $$\sup_{[0, T)\times (-\infty,\infty)} H \leq  \sup_{\{0\}\times (-\infty,\infty)} H \leq C.$$

\qed
\end{proof}

We then have the following immediate corollary.
\begin{proposition} There exists $C>0$ such that on $[0, T) \times (-\infty, \infty)$,

\begin{equation}
 u'' \leq C \min ( T-t, e^{ \rho}, e^{-\rho}).
 \end{equation}

\end{proposition}

\begin{proof}
It suffices to prove that $e^\rho u''$ is bounded above by the previous lemma.
Let $H_\gamma = e^{-t} e^{\gamma \rho} u''$ for $\gamma\in (0,1)$. Then the evolution of $H_\gamma$ is given by

$$\ddt{}H_\gamma = e^{-t} e^{\gamma \rho} (    \frac{u^{(4)}} {u''} - \frac{ (u''')^2 }{ (u'')^2} + \frac{ m u'''}{u'} - \frac{ m (u'')^2}{ (u')^2} + \frac{ n u'''}{ a + u'} - \frac{ n (u'')^2}{( a + u')^2}  ) - H_\gamma.$$
Also for any $t\in [0, T)$, $$\lim_{|\rho|\rightarrow \infty} H_\gamma(\rho) = 0.$$
Suppose $H_\gamma(t_0, \rho_0) = \sup_{[0, t_0]\times (-\infty, \infty)} H_\gamma(t, \rho)$. Then at $(t_0, \rho_0)$, we have $$ u^{(4)}\leq -\gamma u''', ~~~ u'''= -\gamma u''$$
by the maximum principle and then
$$0\leq \ddt{}H_\gamma \leq -H_\gamma.$$
Hence $H_\gamma \leq \sup_{ \rho\in (-\infty, \infty)} H_\gamma(0, \rho)$ and there exists $C>0$ such that for $t\in [0, T)$ and $\gamma\in (0, 1)$, $H_\gamma \leq C.$ By letting $\gamma \rightarrow 1$, we can uniformly bound $e^{-t} e^\rho u''$ on $[0, T) \times (-\infty, \infty)$ from above and the lemma follows.

\qed
\end{proof}

We then obtain uniform bounds for the evolving metrics from the upper bound on $u'$ and $u''$.

\begin{corollary}\label{collapsing metric}
There exists $C>0$ such that on $[0, T)\times X_{m,n}$,
\begin{equation}
a(t)\omega_{FS} \leq \omega(t)\leq (a(t) + C(T-t)) \omega_{FS}+C \min \{ (T-t) (e^{-\rho}+e^\rho) \hat\omega, \hat\omega \}.
\end{equation}

\end{corollary}

\begin{proposition}\label{fiber diameter}
For any $\epsilon>0$, there exists $T_\epsilon\in (0,T)$ such that for $t\in (T_\epsilon, T)$ and any fibre $X_z$ with $z\in \mathbb{P}^n$,

$$diam(X_z, g(t)|_{X_z}) < \epsilon.$$

\end{proposition}

\begin{proof}  We consider the following open set $V_\kappa \subset X_{m,n}$ for $\kappa>0$ defined by
$$ V_\kappa = \{ \kappa^{-1} < e^\rho < \kappa\}.$$
Since $\omega|_{X_z} \leq C  \hat\omega |_{X_z}$, for any $\epsilon>0$, there exists $\kappa_\epsilon>0$ such that for all $t\in [0, T)$ and $\kappa>\kappa_\epsilon$,
$$diam (X_z\cap(X\setminus V_\kappa), g(t)) < \epsilon/2$$
by similar argument in the proof of Theorem \ref{ghcon2}.
On the other hand,  in $V_{2\kappa_\epsilon}$, $\omega |_{X_z} \leq C(T-t)\hat\omega|_{X_z}$. Then there exists $T_\epsilon <T$ such that
$$diam (X_z\cap V_{2\kappa_\epsilon}, g(t)) < \epsilon/2.$$
The proposition then follows easily.

\qed
\end{proof}

\begin{theorem}\label{ghcon} There exists $C>0$ such that for $t\in[0, T)$,
$$diam( X_{m,n}, g(t)) \leq C.$$
Furthermore, $(X_{m,n}, g(t))$ converges to $(\mathbb{P}^n, (a_0- \frac{n-m}{m+2} b_0) \omega_{FS})$ in Gromov-Hausdorff sense as $t\rightarrow T$.

\end{theorem}

\begin{proof}

Let $V_{\kappa} = \{ \kappa^{-1} \leq e^\rho \leq \kappa \}$ for $\kappa >0$. From the calculation above, there exists $C>0$ such that on $V_{\kappa}$,
$$ a\omega_{FS} \leq \omega(t) \leq a\omega_{FS} + C(T-t) \omega_{FS} +  C_\kappa(T-t) \hat\omega  $$
and so $\omega(t)$ converges to $\omega_{FS} $ uniformly in $C^0(V_{\kappa_1} )$ as $t\rightarrow T$. On the other hand, the diameter of any fibre $X_z$ for $z\in \mathbb{P}^n$ tends to $0$ uniformly as $t\rightarrow T$.

We now choose a smooth map $\sigma: \mathbb{P}^m \rightarrow X_{m,n}$ such that the image of $\sigma$ sits in the interior of $V_1$. Then the theorem follows by similar argument in the proof of Theorem 5.1 in \cite{SW1}.

\qed
\end{proof}

\subsection{The case $m < n$ and $\frac{b_0}{m+2} = \frac{a_0}{n-m}$}

In this case, $X_{m,n}$ is Fano and the initial K\"ahler class is proportional to $c_1(X_{m,n})$. The first singular time of the K\"ahler-Ricci flow is $T= \frac{b_0}{m+2} = \frac{a_0}{n-m}$. By Perelman's diameter estimates, we have $$diam(X_{m,n}, g(t)) \leq C(T-t) $$ for a constant $C>0$ and so the flow becomes extinct at $t=T$.

\section{Resolution of singularities by the K\"ahler-Ricci flow}

\subsection{Resolution by the Fubini-Study metric and its Ricci curvature}

Consider the morphism $\Phi_{m,n} : X_{m,n} \rightarrow \mathbb{P}^{(m+1)(n+1)}$ by $\mathcal{O}(1)$.  We assume that $m$, $n\geq 1$.
The restriction of the Fubini-Study metric on $Y_{m,n}$,  the image of $\Phi_{m,n}$, is given by
\begin{eqnarray*}
\hat\omega &=& \ddbar \log (1+ e^\rho) \\
&=& \frac{e^\rho}{1+e^\rho} \omega_{FS} + h_\xi  (\frac{1}{1+e^\rho} \delta_{\alpha \beta}  -   \frac{1}{(1+e^\rho)^2}  h_\xi\xi^{\bar \alpha} \xi^{\beta} )  \frac{\sqrt{-1}}{2\pi}  \nabla \xi^\alpha \wedge \nabla \xi^{\bar \beta}.
\end{eqnarray*}
Its induced volume from on $Y_{m,n}$ is given by
$$\hat\omega^{m+n+1} =    h_\xi ^{m+1}  \frac{e^{n \rho}} {(1+e^\rho)^{m+n+2}  }(\omega_{FS} ^ n \wedge \prod_{\alpha=1}^{m+1}  \frac{\sqrt{-1}}{2\pi}  d\xi^\alpha \wedge d\xi^{\bar\alpha}).  $$
We can now calculate the Ricci form.
\begin{eqnarray*}
&&-Ric(\hat\omega) \\
&=& \ddbar ( n\rho - (m+n+2) \log (1+e^\rho) ) +  (m-n) \omega_{FS} \\
&=& (m-n)\omega_{FS} + \ddbar \hat u_{Ric}\\
&=&(m - \frac{(m+n+2) e^\rho}{1+ e^\rho} )  \omega_{FS} \\%
&&+ \frac{\sqrt{-1}}{2\pi} h_\xi e^{-\rho}\{  \left( n- \frac{  (m+n+2)e^\rho}{1+e^\rho} \right)\delta_{\alpha \beta}  + e^{-\rho} \left( \frac{(m+n+2) e^{2\rho} } {(1+e^\rho)^2} - n  \right) h_\xi \xi^{\bar \alpha} \xi^{\beta} ) \}   \frac{\sqrt{-1}}{2\pi}  \nabla \xi^\alpha \wedge \nabla \xi^{\bar \beta},
\end{eqnarray*}
where $\hat u_{Ric} = n\rho - (m+n+2) \log (1+e^\rho)$.

Let $O$ be the vertex of $Y_{m,n}$ as the projective cone over $\mathbb{P}^m\times \mathbb{P}^n$. We define $(Y_{m,n})_{reg} = Y_{m,n}\setminus \{O \}$ as the nonsingular part of $Y_{m,n}$.
\begin{lemma}
We define $\hat\omega_{\epsilon} = \hat\omega - \epsilon Ric(\hat\omega)$. Then there exists $\epsilon_0>0$, such that $\hat\omega_\epsilon>0$ on $(Y_{m,n})_{reg}$, the nonsingular part of of $Y_{m,n}$ for $\epsilon \in (0, \epsilon_0)$.

\end{lemma}

\begin{proof} It suffices to check for $\rho\leq 0$ because $\hat\omega$ is K\"ahler on $Y_{m,n}\setminus \{ O\}$ and $Ric(\hat\omega)$ is smooth away from $O$. The calculation is straightforward by assuming $\xi= (|\xi|, 0, ..., 0)$ after certain $U(m+1)$ transformation.

\qed
\end{proof}

Let ${\tilde  X}_{m,n}$  be the blow-up of $X_{m,n}$ along the zero section $P_0$. Then we have the following commutative diagram from section 1.9 in \cite{D}.

\begin{equation}
\begin{diagram}
\node{X_{m,n} } \arrow{s,l}{{\small \pi_{m,n}} }     \node{\tilde X_{m,n}}  \arrow{w,t}{\small \vartheta_1} \arrow{s,r}{{\small \Psi}}  \arrow{e,t}{\small \vartheta_2}    \node{ X_{n,m}}  \arrow{s,r}{{\small \pi_{n,m}}}\\
\node{\mathbb{P}^n}      \node{\mathbb{P}^m\times \mathbb{P}^n} \arrow{w,t}{ p_1}   \arrow{e,t}{ p_2}\node{\mathbb{P}^m}
\end{diagram}
\end{equation}

\begin{proposition}\label{ricreso}

The metric completion of $((Y_{m,n})_{reg}, \hat \omega_\epsilon )$ is a smooth compact metric space isomorphic to $\tilde X_{m,n}$ for sufficiently small $\epsilon>0$. In particular, $\hat\omega_\epsilon$ extends to a smooth K\"ahler metric on $\tilde X_{m,n}$.

\end{proposition}

\begin{proof} The potential $\hat u_\epsilon$ of $\hat\omega_\epsilon$ is given by
$$\hat u_\epsilon= \hat u + \epsilon \hat u_{Ric} = n\epsilon \rho  + (1-(m+n+2)\epsilon )\log (1+\rho).$$
It suffices to compare $\hat\omega_\epsilon$ to a smooth K\"ahler metric on $\tilde X_{m,n}$. We let $$\tilde u = a\rho +  b \log (1 + e^\rho)$$ with $a, b>0$. In particular, $\tilde u = \hat u_\epsilon$ when $a=n\epsilon$ and $b= 1-(m+n+2)\epsilon$.
Then $$\tilde\omega = \ddbar \tilde u = ( a+b\frac{e^\rho}{1+e^\rho}) \omega_{FS} + \frac{\sqrt{-1}}{2\pi} e^{-\rho}h_{\xi} ( \tilde u' \delta_{\alpha \beta} + h_\xi e^{-\rho} ( \tilde u'' - \tilde u')  \xi^{\bar \alpha} \xi^{\beta} )\nabla \xi^\alpha \wedge \nabla \xi ^{\bar \beta}.$$
$\tilde \omega$ restricted on each fibre $\mathbb{P}^{n+1}\cap (X_{m,n}\setminus P_0)$ is give by
$$h_\xi e^{-\rho} ( \tilde u' \delta_{\alpha \beta} + e^{-\rho} ( \tilde u'' - \tilde u')  h_\xi \xi^{\bar \alpha} \xi^{\beta} )  \frac{\sqrt{-1}}{2\pi}  \nabla \xi^\alpha \wedge \nabla \xi ^{\bar \beta}$$
whose metric completion is exactly $\mathbb{P}^{n+1}$ blown up at one point. Therefore $\ddbar \tilde u$ blows up along the zero section $P_0$ and replaces $P_0$ by the $\mathbb{P}^1$-bundle over $P_0$ (isomorphic to $\Phi_{m,n}^* \mathbb{P}^m$ blown up at one point) and $\Psi^* (\ddbar \rho)$ is positive on the exceptional divisor $E=\mathbb{P}^m\times \mathbb{P}^n$.

$\ddbar \rho= \ddbar\log ((1+|z|^2) |\xi|^2)$ is exactly the pullback of the Fubini-Study metric on $\mathbb{P}^m\times \mathbb{P}^n$ by $\Psi_{m,n}$.   On the other side, $\tilde\omega$ lies in a K\"ahler class $\Psi^*[\ddbar \rho] + b \Phi_{m,n}^* [ \hat\omega]$ on $\tilde{X}_{m,n}$. Also $\hat\omega$ is positive on $\tilde X_{m,n}\setminus (\mathbb{P}^m\times\mathbb{P}^n)$. Therefore $a\Psi^* \ddbar \rho + b \vartheta^*\hat\omega$ defines a smooth K\"ahler metric on $\tilde{X}_{m,n}$ for $a$, $b>0$ and the proposition follows.

\qed
\end{proof}

 For a given projective embedding  $X \hookrightarrow \mathbb{P}^N$ of a normal variety $X$, the Ricci curvature is well-defined on $X_{reg}$, the nonsingular part of $X$, for the restriction of the Fubini-Study metric $\omega_{FS}$.  We consider $\omega_\epsilon= \omega_{FS}- \epsilon Ric(\omega_{FS})$ for sufficiently small $\epsilon>0$.  Let $(\tilde X, d_\epsilon)$ be the metric completion of $(X_{reg}, \omega_\epsilon)$. Then Proposition \ref{ricreso} suggests  that $\tilde X$ is possibly a resolution of singularities for $X$. However, such a resolution is not necessarily minimal as shown in the example above.  This leads us to consider the K\"ahler-Ricci flow as a certain smoothing process to resolve the singularity of a general normal variety. The goal of the section is to show that indeed the K\"ahler-Ricci flow gives an optimal resolution of singularities for $Y_{m,n}$.

\subsection{The case  $m\neq n$}

In the section, we will consider the K\"ahler-Ricci flow on $Y_{m,n}$ with the initial metric $\omega_0=b_0\hat \omega$ for some $b_0>0$. Since $Y_{m,n}=Y_{n,m}$, we can assume that $m> n$.

 We choose the potential for $\hat \omega$ to be
$$ \hat u= \log (1+ e^\rho). $$
Then calculation in section \ref{redric} suggests that the K\"ahler-Ricci flow should be equivalent to a parabolic PDE for $u$ as below,

\begin{equation}
\ddt{u} = \log [ (a + u')^m (u')^n u''] - (n+1)\rho, ~~~u|_{t=0} =b_0 \hat u
\end{equation}
with $a(t) = (n-m) t$,
or
\begin{equation}\label{correctflow}
\ddt{u} = \log [ (a + u')^n (u')^m u''] - (m+1)\rho, ~~~u|_{t=0} = b_0 \hat u,
\end{equation}
with $a(t)=(m-n)t$, since $a_0=0$.

We have to choose (\ref{correctflow}) because $m> n$ and $a(t)$ should be nonnegative for $t>0$. This can be seen by the class evolution of the K\"ahler-Ricci flow because  $X_{m,n}$ is the only resolution of $Y_{m,n}$ such that $K_{X_{m,n}}$ is $\mathbb{Q}$-Cartier and the class $[\hat \omega] + \epsilon [K_{X_{m,n}}] >0$ for sufficiently small $\epsilon>0$. Hence now we can lift the K\"ahler-Ricci flow on $Y_{m,n}$ to the one on $X_{m,n}$ starting with $\hat\omega$.

Note that $\hat\omega$ has bounded local potential and for any smooth K\"ahler metric $\omega_0$ on $X_{m,n}$, there exists $C>0$ such that $$\hat\omega\leq C \omega_0.$$ By \cite{SoT3}, $\hat\omega\in \mathcal{K}_{[\hat\omega], \infty} (X_{m,n})$ (cf. \cite{SW2}) and there exists a unique weak K\"ahler-Ricci flow on $X_{m,n}$ starting with $\hat\omega$. Furthermore, the solution becomes smooth K\"ahler metrics on $X_{m,n}$ once $t>0$. Therefore, it suffices to study the behavior of the solution as $t\rightarrow 0^+$.

We first write down the equivalent parabolic flow of Monge-Ampere type for the K\"ahler-Ricci flow. Since $[\hat\omega]+ t [K_{X_{m,n}}]>0$ for sufficiently small $t>0$, there exists a smooth volume form $\Omega$ with $\chi=\ddbar\log \Omega$, such that  $$\omega_t = \hat\omega + t\chi>0,$$
for $t\in (0, T)$, where $T=\sup\{t>0~|~[\hat\omega]+ t [K_{X_{m,n}}]>0\}. $ Let the solution of the K\"ahler-Ricci flow be given as $\omega(t) = \omega_t + \ddbar\varphi$. Then

\begin{equation}\label{singfl}
\ddt{\varphi} = \log \frac{ (\omega_t + \ddbar \varphi)^{m+n+1}}{\Omega}, ~~~\varphi|_{t=0} =0.
\end{equation}

It is proved in \cite{SoT3}, that $\varphi\in C^\infty( (0, T)\times X_{m,n})\cap C^\infty([0, T)\times (X_{m,n}\setminus P_0))$ and

$$||\varphi||_{L^\infty([0, T/2]\times X_{m,n})} < \infty.$$

\begin{lemma} \label{volest3} Then there exists $C>0$, such that on $[0, T/2]\times X_{m,n}$,
\begin{equation}
\omega^{m+n+1} \leq C \max\{1, e^{ -(m-n)\rho}\} \Omega,
\end{equation}
and on $[T/2, T)\times X_{m,n}$,
$$\omega^{m+n+1} \leq C \Omega.$$
\end{lemma}

\begin{proof} It suffices to show that the lemma holds on $[0, T/2]\times (-\infty, 0]$, as one can easily obtain the estimate $\omega^{m+n+1} \leq C_1 \Omega$ away from the zero section $P_0$ (see [ST3]), as well as for $t \geq T/2$ (see [SW1]), for $C_1>0$.

Let $v= \ddt{u}$. Then the evolution of $v$ is given by

\begin{equation}
\ddt{v} = \frac{nv'}{a+u'}+\frac{mv'}{u'}+ \frac{v''}{u''} +\frac{n(m-n)}{a+u'}.
\end{equation}

Let $H= e^{-t} (v + (m-n)\rho).$ Then by (\ref{correctflow}), $H(0) \leq C_2 + m\rho \leq C_2$ on $\rho \in (-\infty, 0]$ for $C_2>0$. One can calculate the evolution for $H$,
\begin{equation}
\ddt{H} = \frac{nH'}{a+u'}+\frac{mH'}{u'}+ \frac{H''}{u''}-H - \frac{m(m-n)e^{-t}}{ u'} \leq \frac{nH'}{a+u'}+\frac{mH'}{u'}+ \frac{H''}{u''}  -H.
\end{equation}
%
One notices that $\lim_{\rho \rightarrow -\infty} H(t, \rho) = -\infty$ for any $t$ by Proposition \ref{calabi}. Hence the maximum of $H(t, \cdot)$ is achieved away from $P_0$ for $t \geq 0$. Since the K\"{a}her-Ricci flow is smoothly defined away from $P_0$, it follows from the maximum principle that $H \leq C_3$ on $[0, T/2]\times (-\infty, 0]$ for $C_3>0$. Therefore by (\ref{correctflow}) again, there exists $ C_4>0$, such that

 $$(a+u')^n(u')^m u''=e^{v + (m+1)\rho} = e^{e^t H+(n+1)\rho} \leq C_4 e^{(n+1)\rho}.$$
On the other hand, there exists $C_5>0$ such that on $(-\infty, 0]$:

$$(1+\hat{u}')^n(\hat{u}')^m \hat{u}'' \geq C_5 e^{(m+1)\rho}.$$
Combining them, there exists $C_6>0$, such that on $(-\infty, 0]$,
$$\omega^{m+n+1} \leq C e^{ -(m-n)\rho} \Omega.$$
\qed
\end{proof}

\begin{lemma}\label{schwarz3} There exists $C>0$ such that on $[0, T/2]\times X_{m,n}$,

\begin{equation}
\omega \geq C \hat \omega.
\end{equation}

\end{lemma}

\begin{proof} Let $\theta$ be a smooth K\"ahler metric on $X_{m,n}$. We consider the K\"ahler-Ricci flow on $X_{m,n}$ with initial metric $$\omega_{\epsilon, 0}|_{t=0}= b_0 \hat\omega + \epsilon \theta.$$
Then by same argument as in \cite{SW1, So}, we can show that there exists $C>0$ such that for any $\epsilon\in (0, 1)$, the solution $\omega_\epsilon(t)$ of the K\"ahler-Ricci flow is bounded below by
$$\omega \geq C \hat\omega$$
on $t\in (0, T/2]\times X_{m,n}$.

\qed\end{proof}

\begin{proposition} There exist $A_1$ and $A_2>0$ such that for $t\in [0, T/2]$ and $\rho\leq 0$,

$$ u'' \leq  A_1 u' \leq A_2 ~e^{\frac{n+1}{m+n+1} \rho},$$
%

\end{proposition}

\begin{proof}

By the volume comparison in Lemma \ref{volest3}, there exist $C_1, C_2>0$ such that for $t\in [0, T/2)$ and $\rho \leq 0$,

$$ (a+ u')^m (u')^{n} u'' \leq C_1 e^{-(m-n)\rho} (1+ \hat u')^n (\hat u')^{m} \hat u'' \leq C_2 e^{(n+1)\rho}.$$
Applying the same argument in Corollary \ref{1stder1}, there exists $C_3>0$ such that for $t\in [0, T)$ and $\rho\leq 0$,

$$ u' \leq C_3 e^{\frac{n+1}{m+n+1} \rho}.$$
%

Let $H= \log u'' - \log u'$. To prove $H$ is uniformly bounded from above, one just imitates the argument as in Proposition \ref{2ndder1} and in addition checks at $t=0$ when $u = b_0 \hat u $,

$$H(0) = \log \hat u'' - \log \hat u'= - \log (1+ e^\rho) \leq 0.$$

%

\qed
\end{proof}

\begin{corollary} There exists $C>0$ such that on $[0, T/2]\times X_{m,n}$,

\begin{equation}
C^{-1}  \hat \omega \leq \omega \leq C\{ a \omega_{FS} + \hat\omega + e^{-\frac{m}{m+n+1}\rho } \hat\omega  \},
\end{equation}
where $a= (m-n)t.$
\end{corollary}

\begin{theorem}\label{ghcon3}  $(X_{m,n}, g(t))$ converges to $(Y_{m,n}, \hat g)$ in Gromov-Hausdorff sense as $t\rightarrow 0^+$.

\end{theorem}

\begin{proof}

It is proved in \cite{SoT3} that $g(t)$ converges to $\hat g$ in $C^\infty$ topology of $X_{m,n}\setminus P_0$.
 Let $U_\kappa = \{ e^\rho \leq \kappa\}$ be the $\kappa$-tubular neighborhood of $P_0$. Then it suffices to show that for any $\epsilon>0$, there exist $\kappa_\epsilon>0$ and $T_\epsilon\in (0, T/2]$ such that for any $\kappa<\kappa_\epsilon$ and $t\in (0, T_\epsilon)$, $$diam( U_\kappa\setminus P_0, g(t)) < \epsilon.$$ This can be proved by similar argument as in the proof of Theorem \ref{ghcon2}.

\qed
\end{proof}

\begin{theorem}$(X_{m,n}, g(t))$ converges to $(\mathbb{P}^n, \frac{(m-n)b_0}{m+2}\omega_{FS})$ in  Gromov-Hausdorff sense as $t\rightarrow T^-$.
\end{theorem}

\begin{proof}It follows directly from the equation (\ref{correctflow}) that the singular time $T=\frac{b_0}{m+2}$ and $a(T)=\frac{(m-n)b_0}{m+2}$. The theorem thus follows from Theorem \ref{ghcon} as the K\"{a}hler-Ricci flow on $Y_{m,n}$ becomes the K\"{a}hler-Ricci flow on $X_{m,n}$ after arbitrary short time $t>0$.
\qed
\end{proof}

\subsection{ The case  $m=n$}

We now consider the K\"ahler-Ricci flow on $Y_{n,n}$ starting with $b_0\hat\omega= b_0 \ddbar \hat u$, where $\hat u = \log (1+ e^\rho)$. We would like to lift the flow to the one on $X_{n,n}$. If we let $$T=\sup\{t>0~|~ [\hat\omega]+t[K_{X_{n,n}}]~\textnormal{is~big~and~semi-ample}\},$$
then $$T=\frac{b_0}{n+2}>0.$$

By \cite{SoT3}, the K\"ahler-Ricci flow can always be lifted to the one on $X_{n,n}$ for $t\in [0, T)$. However the solution is in general not smooth since $b_0 [\hat\omega]+t[K_{X_{n,n}}]=(b_0- (n+2) t)[\hat\omega]$ vanishes on $P_0$ of $X_{n,n}$ for any $t\geq0$.

We apply the same method in \cite{SoT3} by approximating the flow (\ref{singfl}) by the smooth data.  We consider the family of flows for $\delta\in (0,1)$,
\begin{equation}\label{appfl}
\ddt{\varphi_\delta} = \log \frac{ (\omega_t + \delta\omega_{FS}+ \ddbar\varphi_\delta)^{m+n+1}}{\Omega}, ~~~ \varphi_\delta |_{t=0}=0,
\end{equation}
where $\omega_{FS}$ is the pullback of the Fubini-Study metric on $\mathbb{P}^n$, $\omega_t= (b_0-(n+2))\hat\omega$ and $\Omega$ is a smooth volume form on $X_{n,n}$ with $\ddbar\log \Omega= -(n+2)\hat\omega$.

The above perturbed flow is equivalent to the following family of parabolic flows,

\begin{equation}
\ddt{u_\delta}= \log[(\delta + u_\delta')^n(u_\delta')^n u_\delta''] - (n+1)\rho, ~~u_\delta |_{t=0}=b_0 \hat u.
\end{equation}

\begin{lemma}  There exists $C>0$ such that for $t\in [0, T)$ and   $\delta\in (0, 1)$,

$$ (\delta + u_\delta')^n (u_\delta')^n u_\delta'' \leq C (1+ \hat u')^n (\hat u')^n \hat u''.$$

\end{lemma}

\begin{proof} It suffices to prove the lemma for $\rho \leq 0$ as the volume estimate holds true away from the zero section $P_0$(see [ST3]). Let $v_\delta = \ddt{u_\delta}$. Then

$$ \ddt{v_\delta} = \frac{n v_\delta'}{\delta + u_\delta'}+ \frac{n v_\delta'}{u_\delta'}+ \frac{v_\delta''}{u_\delta''}  .$$
Let $H_{\delta, \epsilon} = e^{-t} ( v_\delta +\epsilon \rho)$ for $\epsilon\in (0,1)$. Then there exists $C_1>0$ such that $\lim_{\rho\rightarrow -\infty} H_{\delta,\epsilon} \leq C_1$ for fixed $\delta$, $\epsilon$ and $t\in [0, T)$.
$$ \ddt{H_{\delta, \epsilon}} = \frac{n H_{\delta, \epsilon}'}{\delta + u_{\delta}'}+ \frac{n H_{\delta, \epsilon}'}{u_\delta'}+ \frac{H_{\delta, \epsilon}''}{u_\delta''}  - \frac{n \epsilon e^{-t}}{\delta + u_{\delta}'} - \frac{n \epsilon e^{-t}}{u_\delta'}- H_{\delta, \epsilon}.$$
The maximum principle implies that there exists $C_2>0$ such that for $t\in[0, T)$, $\delta\in (0,1)$, and $\epsilon\in (0,1)$, $$\sup_{t\in [0,T), \rho\in (-\infty,0]} H_{\delta, \epsilon} \leq  \sup_{t=0, \rho\in (-\infty, 0]} H_{\delta, \epsilon} + C_2.$$
The lemma is then proved by checking $ H_{\delta, \epsilon}(0, \cdot)$ is bounded from above and letting $\epsilon\rightarrow 0$.

\qed
\end{proof}

The following proposition can be proved in the same way as in Proposition \ref{2ndder1}.
\begin{proposition} There exist $C_1, C_2>0$ such that for $t\in [0, T)$, $\rho\in(-\infty, \infty)$ and  $\delta\in (0, 1)$,
\begin{equation}\label{appr}  u_\delta'' \leq C_1 u_\delta' \leq C_2 \min(1, e^{\frac{n+1}{2n+1} \rho} ).
\end{equation}

\end{proposition}

\begin{corollary} There exists $C>0$ such that for $t\in [0, T)$, $\delta\in (0, 1)$ and $\rho\leq 0$,

\begin{equation}\label{vol4} \omega (t) \leq \delta \omega_{FS} + C e^{-\frac{n}{2n+1} \rho} \hat \omega . \end{equation}
Furthermore,  for any $t\in [0, T)$, there exists $C_t>0$ such that \begin{equation}\label{schwarz4}\omega(t) \geq  C_t \hat \omega.\end{equation}

\end{corollary}

\begin{proof} Letting $\delta\rightarrow 0$ in equation (\ref{appr}), we have $u'' \leq C_1 u' \leq C_2 e^{\frac{n+1}{2n+1} \rho}$, for constants $C_1, C_2>0$. Equation (\ref{vol4}) then follows easily.

For any $T'\in(0,T)$, we can apply the argument in Lemma \ref{schwarz3} to show that there exists $C_{T'}>0$ such that for $\delta\in(0,1)$ and on $[0, T')\times X_{n,n}$, $$\omega_\delta=\delta \omega_{FS}+\ddbar u_\delta \geq C_{T'} \hat \omega.$$ Inequality (\ref{schwarz4}) follows by letting $\delta \rightarrow 0$.
\qed
\end{proof}

Let $(X_t, d_t)$ be the metric completion of $(X_{n,n}\setminus P_0, \omega(t))$ for $t\in (0, T)$.

\begin{theorem} For any $t\in (0, T)$, $(X_t, d_t)$ has finite diameter and $(X_t, d_t)$ is homeomorphic to the projective cone $Y_{n,n}$ over $\mathbb{P}^n\times \mathbb{P}^n$ via the Segre map. Furthermore, the Gromov-Hausdorff distance $D(t)=d_{GH}( (X_t, d_t), (Y_{n,n}, \hat g))$ is a  continuous function in $t\in [0, T)$ and
\begin{equation}
\lim_{t\rightarrow 0} D(t)=0.
\end{equation}

\end{theorem}

\begin{proof} We will consider the approximating K\"ahler-Ricci flow defined by (\ref{appfl}). The solution $\omega_\delta(t)=\omega_t + \delta \theta+\ddbar \varphi_\delta$ is a smooth K\"ahler metric on $(0, T)\times X_{n,n}$. Let $U_\kappa = \{ e^\rho \leq \kappa\}$ be the $\kappa$-tubular neighborhood of $P_0$. By similar argument in the proof of Theorem \ref{ghcon2}, we can show that for any fixed $t\in (0, T)$ and $\epsilon>0$, there exist $\kappa_\epsilon>0$ and $\delta_\epsilon>0$, such that for any $\kappa\in (0, \kappa_\epsilon) $ and $\delta\in (0, \delta_\epsilon]$,
$$diam( U_{\kappa}, g_\delta(t))< \epsilon.$$
Since $g_\delta(t)$ converges to $g(t)$ in $C^\infty$-topology on $(0, T) \times X_{n,n}\setminus P_0$, we can show by similar argument in \cite{SW1} that $(X_t, d_T)$ has finite diameter and is homeomorphic to $(Y_{n,n}, \hat g)$.

Now it suffices to show that  $D(t)$ is continuous and the rest of the theorem can be proved by similar argument in the proofs of Theorem \ref{ghcon2} and Theorem \ref{ghcon3}. Fix $t_0\in (0, T)$, we consider
$$\hat D(t) = d_{GH} ((X_t, d_t), (X_{t_0}, d_{t_0}))$$ for $t\in [0, T)$. We claim that $$\lim_{t\rightarrow t_0} \hat D(t) =0 .$$
First note that $\omega(t)$ converges to $\omega(t_0)$ in $C^\infty$ topology of $X_{n,n}\setminus P_0$ as $t\rightarrow t_0$. On the other hand, we can apply the same argument as before that, for any $\epsilon>0$, there exists $\eta>0$ and $\kappa_\epsilon>0$ such that for any $\kappa< \kappa_\epsilon$ and $t\in (t_0-\eta, t_0+\eta) \cap (0, T/2]$,
$$ diam( U_\kappa\setminus P_0, g(t)) < \epsilon.$$

Also for any $\epsilon>0$, any $\kappa>0$ and $k>0$ , there exists $\eta>0$  such that for $t\in [t_0-\eta, t_0+\eta) \cap (0, T/2]$,

$$|| \omega(t)- \omega(t_0)||_{C^k(X_{n,n}\setminus U_\kappa)}< \epsilon.$$ Here the $C^k$ norm is taken with respect to a fixed K\"ahler metric $\theta$ on $X_{n,n}$. Then the claim follows by similar argument in \cite{SW1}.

\qed
\end{proof}

\subsection{Finite time extinction on $Y_{n,n}$}

In this section, we consider the limiting behavior of the K\"ahler-Ricci flow on $Y_{n,n}$ starting with $b_0\hat \omega$ for some $b_0>0$. We have shown the existence of the solution $g(t)$ as in section \ref{col}.

\begin{theorem} Let $(Y_{n,n}, d_t)$ be the metric completion of $(Y_{n,n}\setminus \{O\}, g(t))$. Then

$$\lim_{t\rightarrow T} diam( Y_{n,n}, d_t) =0.$$

\end{theorem}

\begin{proof} We again consider the perturbed flow (\ref{appfl}). Notice that the K\"ahler class along the flow is given by $ \delta[\omega_{FS}]+(b_0-(n+2)t) [\hat\omega]$, hence for all $\delta\in(0,1)$, we have
$$0<u_\delta'\leq (b_0-(n+2))t= (n+2)(T-t).$$

By similar argument as in the section \ref{collapse}, there exist $C_1$ and $C_2>0$ such that on $[0, T)\times(-\infty, \infty)$
$$u_\delta''\leq C_1 u_\delta' \leq C_2 \min(T-t, e^{\frac{n+1}{2n+1}\rho}), ~~~~u_\delta'' \leq C_2 e^{-\rho}$$
for $\delta\in(0,1)$. By letting $\delta\rightarrow 0$, we have
$$u''\leq C_1 u' \leq C_2 \min(T-t, e^{\frac{n+1}{2n+1}\rho}), ~~~~ u''\leq C_2 e^{-\rho}.$$
We thus have the estimates on $Y_{n,n}\setminus \{ O\}$,

$$\omega(t) \leq C  \min ( (T-t)(e^{-\rho} + e^\rho) \hat\omega  , (1+ e^{-\frac{n}{2n+1}\rho}) \hat\omega ).$$
Then by similar argument in \cite{SW1}, we can show that
$$\lim_{t\rightarrow T} diam(Y_{n,n}\setminus \{O\}, \omega(t)) =0.$$ The theorem follows since $(Y_{n,n}, d_t)$ is the metric completion of $( Y_{n,n}\setminus \{ O \}, g(t) )$.

\qed
\end{proof}


\bigskip

\noindent {\bf{Acknowledgements:}} The authors are grateful to Professor
D. H. Phong for his support and encouragement. They also thank Professor G. Tian, Chenyang Xu and Chi Li for many helpful discussions. In addition, the authors would like to thank V. Tosatti for a number of useful suggestions. The second named author is also grateful to Professor X. Huang for his constant support. Part of the work was carried out during the first-named author's visit at Columbia University, and he thanks the department for its kind hospitality. It was also carried out during the second named author's graduate study at  Rutgers University.


\bigskip
\bigskip

\noindent $^{*}$ Department of Mathematics, Rutgers University, New Brunswick, NJ 08854\\

\noindent $^\dagger$ Department of Mathematics, Johns Hopkins University, Baltimore, MD 21218\\

\end{document}